\documentclass[reqno]{amsart}
\usepackage{amsmath,amssymb,amsthm,ascmac}
\usepackage[mathscr,mathcal]{eucal}
\usepackage{latexsym}
\usepackage{enumerate}
\usepackage[dvipdfmx]{graphicx}

\begin{document}

\allowdisplaybreaks

\theoremstyle{definition}
\newtheorem{defi}{\textbf{Definition}}[section]
\newtheorem{thm}[defi]{\textbf{Theorem}}
\newtheorem{lem}[defi]{\textbf{Lemma}}
\newtheorem{prop}[defi]{\textbf{Proposition}}
\newtheorem{cor}[defi]{\textbf{Corollary}}
\newtheorem{ex}[defi]{\textbf{Example}}
\newtheorem{rem}[defi]{\textbf{Remark}}
\newtheorem*{corr}{\textbf{Corollary}}

\theoremstyle{plain}
\newtheorem{maintheorem}{Theorem}
\newtheorem{theorem}{Theorem }[section]
\newtheorem{proposition}[theorem]{Proposition}
\newtheorem{mainproposition}{Proposition}
\newtheorem{lemma}[theorem]{Lemma}
\newtheorem{corollary}[theorem]{Corollary}
\newtheorem{maincorollary}{Corollary}
\newtheorem{claim}{Claim}
\renewcommand{\themaintheorem}{\Alph{maintheorem}}
\theoremstyle{definition} \theoremstyle{remark}
\newtheorem{remark}[theorem]{Remark}
\newtheorem{example}[theorem]{Example}
\newtheorem{definition}[theorem]{Definition}
\newtheorem{problem}{Problem}
\newtheorem{question}{Question}
\newtheorem{exercise}{Exercise}

\newtheorem*{subject}{�ړI}
\newtheorem*{mondai}{Problem}
\newtheorem{lastpf}{�ؖ�}
\newtheorem{lastpf1}{proof of theorem2.4}
\renewcommand{\thelastpf}{}

\newtheorem{last}{Theorem}

\renewcommand{\thelast}{}
\renewcommand{\proofname}{\textup{Proof.}}

\renewcommand{\theequation}{\arabic{section}.\arabic{equation}}
%%%%%  ���ԍ��� section ���Ƀ��Z�b�g����  %%%%%
\makeatletter
\@addtoreset{equation}{section}

\title[Topological entropy of generic sets for ($\alpha$-$\beta$)-shifts]
{Topological entropy of the set of generic points for ($\alpha$-$\beta$)-shifts}
\author[K. Yamamoto]{Kenichiro Yamamoto}
\address{Department of General Education \\
Nagaoka University of Technology \\
Nagaoka 940-2188, JAPAN
}
\email{k\_yamamoto@vos.nagaokaut.ac.jp}


\date{}

\maketitle
\large

\begin{abstract}
We prove that all ($\alpha$-$\beta$)-shifts with $0\le \alpha<1$ and $\beta>2$ are saturated, that is,
for any invariant measure, the topological entropy of the set of generic points coincides with the metric entropy.
\end{abstract}

\section{Introduction}

Let $(X,d)$ be a compact metric space
and $f:X\to X$ be a continuous map.
For an $f$-invariant Borel probability measure $\mu$,
we denote by $G_{\mu}$ the set of points
whose asymptotic distributions are equal to $\mu$, i.e.
$$G_{\mu}:=\left\{
x\in X ~: 
\lim_{n\to \infty} \frac{1}{n} \sum_{j=0}^{n-1}
\delta_{f^j(x)}= \mu  \right\}.$$
Here $\delta_y$ stands for the $\delta$-measure at the point $y\in X$.
Then, it is well-known that
$\mu (G_{\mu})=1$ if ergodic and $\mu(G_{\mu})=0$ otherwise.
The points in $G_{\mu}$ are called {\it generic points} of $\mu$.

In multifractal analysis, 
the 'size' of $G_{\mu}$ has been estimated 
by several means, 
including the Hausdorff dimension (\cite{HS,PS3}), the entropy (\cite{FLP,HTW,PS2,TV}) and
the topological pressure (\cite{Y,ZC}).
In this paper we consider the topological entropy of $G_{\mu}$.
The notion of topological entropy for non-compact sets 
is introduced by Bowen in \cite{B}.
By using this definition he proved that if $\mu$ is ergodic, then
\begin{equation}
\label{eq1.1}
h_{top}(f,G_\mu)=h(\mu)
\end{equation}
where $h_{top}(f,G_{\mu})$ is the topological entropy of $G_{\mu}$
and $h(\mu)$ is the metric entropy of $\mu$. 
However, when $\mu$ is non-ergodic,
we can easily construct an example which does not satisfy the equation (\ref{eq1.1}).

The system $(X,f)$ is said to be \textit{saturated} if (\ref{eq1.1}) holds
for any $f$-invariant measure $\mu$ on $X$.

Recently, several authors investigated the saturatedness for important classes of dynamical systems.
In \cite{FLP}, Fan, Liao and Peyrire proved several dynamical properties implied by the saturatedness,
and also show the saturatedness for systems with the specification property.
At the same time, Pfister and Sullivan in \cite{PS2} prove the saturatedness
for more general systems beyond specification, including all $\beta$-shifts.
Motivated by above results,
in this paper, we investigate the saturatedness for ($\alpha$-$\beta$)-shifts.

The ($\alpha$-$\beta$)-shift was introduced by Parry in \cite{Pa} as a natural, but
non-trivial generalization of the $\beta$-shift,
and its dynamical properties has been studied by many authors (\cite{FP,FL,H3,H,LSS,RS}).
Let $0\le \alpha<1$ and $\beta>1$. An ($\alpha$-$\beta$)-transformation
$T_{\alpha,\beta}\colon [0,1)\to [0,1)$ is defined by
\begin{equation}
\label{alpha-beta}
T_{\alpha,\beta}(x):=
\beta x+\alpha\ \ (\text{mod}\ 1)
\end{equation}
The ($\alpha$-$\beta$)-shift $\Sigma_{\alpha,\beta}^+$ is a subshift
consisting of all kneading sequences induced by $T_{\alpha,\beta}$ (see \S2.3 for the precise definition).
If $\alpha=0$, then $\Sigma_{\alpha,\beta}^+$ coincides with the $\beta$-shift.
We remark that the ($\alpha$-$\beta$)-transformation is also called
the \textit{intermediate $\beta$-transformation}, or the \textit{linear mod $1$ transformation}.
Similarly, the ($\alpha$-$\beta$)-shift is
also called the \textit{intermediate $\beta$-shift}.

While the definition of  ($\alpha$-$\beta$)-shift is similar to that of $\beta$-shift,
some of dynamical properties are different. For example,
it is known that all $\beta$-shifts are topologically mixing, but ($\alpha$-$\beta$)-shift is not transitive
for some $1<\beta<2$ (see \cite{RS}).
Moreover, it is not difficult to see that the saturatedness is failed for certain non-transitive
($\alpha$-$\beta$)-shifts. Hence it is natural to ask which ($\alpha$-$\beta$)-shifts are saturated.
In this paper, we prove the saturatedness for ($\alpha$-$\beta$)-shifts under the mild condition
``$\beta>2$".
Now, we state our main result of this paper.

\begin{maintheorem}
\label{main}
Let $0\le\alpha<1$, $\beta>2$, $\Sigma_{\alpha,\beta}^+$ be a ($\alpha$-$\beta$)-shift
and $\sigma\colon \Sigma^+_{\alpha,\beta}\to \Sigma^+_{\alpha,\beta}$ be a shift map.
Then $(\Sigma^+_{\alpha,\beta},\sigma)$ is saturated.
\end{maintheorem}

In the proof, we show the entropy density of ergodic measures for ($\alpha$-$\beta$)-shifts
(Proposition \ref{entropy-dense}), which plays a very important role to prove Theorem \ref{main}.
Here the entropy density means that every invariant measure can be
approximated by ergodic measures with similar entropies.
To show the entropy density,
we use the density of periodic measures and
Hofbauer's Markov Diagram.
This is a quite new technique, which was first established in Chung and the present author's recent preprint \cite{CY}
to obtain the large deviation bounds for ($\alpha$-$\beta$)-transformations.

This paper is organized as follows:
In \S2, we give our definitions, prove several basic lemmas and recall some facts of ($\alpha$-$\beta$)-shifts.
We give a proof of Theorem \ref{main} in \S3.

\section{Preliminaries}
\subsection{Topological entropy for non-compact sets}
Let $(X,d)$ be a compact metric space and $f\colon X\to X$ be a continuous map.
For $n\in\mathbb{N}$, $x\in X$ and $\epsilon>0$, we set
$$B_n(x,\epsilon):=\left\{y\in X:\max_{0\le k\le n-1}d(f^k(x),f^k(y))\le\epsilon\right\}.$$
In this subsection, we recall the definition of the topological entropy for non-compact set, which was introduced by Bowen in \cite{B}.
For $s\in\mathbb{R}$, $M\in\mathbb{N}$ and $\epsilon>0$, we set
$$M(Z,s,M,\epsilon):=\inf_{\Gamma}\left\{\sum_{B_{n_i}(x_i,\epsilon)\in\Gamma}e^{-sn_i}\right\}$$
where the infimum is taken over all $\Gamma:=\{B_{n_i}(x_i,\epsilon)\}$ which is a finite or countable cover of $Z$. Since the quantity $M(Z,s,M,\epsilon)$ does not decrease with $M$, we can define
$$M(Z,s,\epsilon):=\lim_{M \rightarrow \infty}
M(Z,s,M,\epsilon).$$
Then we can easily see that the following critical value exists:
$$h_{top}(f,Z,\epsilon):=\inf\{s:M(Z,s,\epsilon)=0\}=\sup\{s:M(Z,s,\epsilon)=\infty\}.$$
Finally we set
$$h_{top}(f,Z):=\lim_{\epsilon \rightarrow 0}h_{top}(f,Z,\epsilon)$$
and call this the \textit{topological entropy} of $Z$.
Then it is easy to see that $h_{\rm top}(f,Z_1)\le h_{\rm top}(f,Z_2)$ if $Z_1\subset Z_2\subset X$.
Moreover, it is proved by Bowen in \cite{B} that the following upper bound for $h_{\rm top}(f,G_{\mu})$
holds:

\begin{lem}(\cite[Theorem 2]{B})
\label{upper}
For any $f$-invariant measure $\mu$ on $X$, we have $h_{\rm top}(f,G_{\mu})\le h(\mu)$.
\end{lem}

\subsection{Symbolic dynamics}
Let $A$ be an at most countable set and
$A^{\mathbb{N}}$ (resp. $A^{\mathbb{Z}}$) be
the set of all one-sided (resp. two-sided) infinite sequences on the alphabet $A$ endowed with
the product topology of the discrete topology on $A$.
The shift map on $A^{\mathbb{N}}$ is $\sigma\colon \omega_1\omega_2\cdots\mapsto
\omega_2\omega_3\cdots$ and the shift map on $A^{\mathbb{Z}}$ is defined analogously.
A one-sided (resp. two-sided) subshift is a closed $\sigma$-invariant set $\Sigma\subset A^{\mathbb{N}}$
(resp. $\Sigma\subset A^{\mathbb{Z}}$).

Let $\Sigma$ be a subshift. The language of $\Sigma$, denoted by $\mathcal{L}(\Sigma)$ is defined by
$$\mathcal{L}(\Sigma):=\left\{w\in \bigcup_{n=1}^{\infty}A^n:[w]\not=\emptyset\right\}.$$
Here $[w]:=\{x\in \Sigma:x_i=w_i,1\le i\le |w|\}$ is a \textit{cylinder set} and $|w|$ is the length of $w$.
We also set $\mathcal{L}_n(\Sigma):=\{w\in\mathcal{L}(\Sigma):|w|=n\}$ for $n\in\mathbb{N}$.
For $u,v\in\mathcal{L}(\Sigma)$, we write
$uv=u_1\ldots u_{|u|}v_1\ldots v_{|v|}$.
We say that a subshift $\Sigma$ is \textit{transitive} if for any $u,v\in\mathcal{L}(\Sigma)$,
there is $w\in\mathcal{L}(\Sigma)$ such that $uwv\in\mathcal{L}(\Sigma)$.

For a subshift $\Sigma$, we denote by $\mathcal{M}(\Sigma)$ the set of all
Borel probability measure on $\Sigma$ endowed with the weak$^{\ast}$-topology.
We also denote by $\mathcal{M}_{\sigma}(\Sigma)\subset\mathcal{M}(\Sigma)$ the
set of all invariant ones, and by $\mathcal{M}_{\sigma}^e(\Sigma)\subset\mathcal{M}_{\sigma}(\Sigma)$,
the set of all ergodic ones. For $\mu\in\mathcal{M}_{\sigma}(\Sigma)$,
we define the \textit{metric entropy} $h(\mu)$ of $\mu$ by
$$\lim_{n\rightarrow\infty}\frac{1}{n}\log\sum_{w\in\mathcal{L}_n(\Sigma)}-\mu[w]\log\mu[w].$$

Let $M=(M_{i,j})_{(i,j)\in A^2}$ be a matrix of zeros and ones.
A (two-sided) \textit{countable Markov shift} $\Sigma_M\subset A^{\mathbb{Z}}$
generated by the transition matrix $M$ is defined by
$$\Sigma_M:=\{x\in A^{\mathbb{Z}}:M_{x_i,x_{i+1}}=1\text{ for }i\in\mathbb{Z}\}.$$ 
Very recently, Takahasi \cite{T} proved the entropy density of ergodic measures for
transitive countable Markov shift as follows:

\begin{lem}(\cite[Main Theorem]{T})
\label{takahasi}
Let $\Sigma_M\subset A^{\mathbb{Z}}$ be a transitive countable Markov shift.
Suppose that $\mu\in\mathcal{M}_{\sigma}(\Sigma_M)$ and $h(\mu)<\infty$.
Then for any $\epsilon>0$ and any neighborhood
$U\subset \mathcal{M}(\Sigma_M)$ of $\mu$, we can find a finite set $F\subset A$ and
an ergodic measure $\rho\in U\cap \mathcal{M}_{\sigma}(F^{\mathbb{Z}}\cap \Sigma_M)$ such that $|h(\mu)-h(\rho)|\le \epsilon$ holds.
\end{lem}

\subsection{($\alpha$-$\beta$)-shifts}
In this subsection, we give the definition of ($\alpha$-$\beta$)-shifts.
In the rest of this paper, we always assume the assumption of Theorem \ref{main}, i.e.
$0\le\alpha <1$ and $\beta>2$.
Let $T_{\alpha,\beta}\colon [0,1)\to [0,1)$ be an ($\alpha$-$\beta$)-transformation
defined by (\ref{alpha-beta}). We set
$$I_1:=\left(0,\frac{1-\alpha}{\beta}\right),I_2:=\left(\frac{1-\alpha}{\beta},\frac{2-\alpha}{\beta}\right),\ldots,
I_k:=\left(\frac{k-1-\alpha}{\beta},1\right),$$
where $k$ denotes the smallest integer, which is not less than $\alpha+\beta$.
We set
$X:=\bigcap_{n=0}^{\infty}T^{-n}\bigcup_{j=1}^kI_j$.
For $x\in X$,
we define a kneading sequence ${\rm It}(x)\in\{1,\ldots,k\}^{\mathbb{N}}$ under $T_{\alpha,\beta}$
by
$${\rm It}(x)_i=j\text{ if }T_{\alpha,\beta}^{i-1}(x)\in I_j.$$
An \textit{($\alpha$-$\beta$)-shift} $\Sigma^+_{\alpha,\beta}$ is then defined by
$$\Sigma^+_{\alpha,\beta}:={\rm cl}(\{{\rm It}(x):x\in X\}).$$
Clearly, $\Sigma^+_{\alpha,\beta}$ is a subshift.
We also define a two-sided ($\alpha$-$\beta$)-shift $\Sigma_{\alpha,\beta}$ by
$$\Sigma_{\alpha,\beta}:=\{x\in\{1,\ldots,k\}^{\mathbb{Z}}:x_ix_{i+1}\cdots\in\Sigma^+_{\alpha,\beta},i\in\mathbb{Z}\}.$$
We say that $\mu\in\mathcal{M}_{\sigma}(\Sigma_{\alpha,\beta})$ is a \textit{periodic measure}
if there exist $x\in\Sigma_{\alpha,\beta}$ and $n\in\mathbb{N}$ such that
$\sigma^n(x)=x$ and $\mu=\frac{1}{n}\sum_{j=0}^{n-1}\delta_{\sigma^j(x)}$.
We denote by $\mathcal{M}_{\sigma}^p(\Sigma_{\alpha,\beta})$ the set of
all periodic measures on $\Sigma_{\alpha,\beta}$.
As we said in the end of \S1, we use the density of periodic measures for ($\alpha$-$\beta$)-shifts,
which is essentially proved by Hofbauer in \cite{H}.

\begin{lem}(\cite[Theorem 2]{H})
\label{per-dense}
The set $\mathcal{M}_{\sigma}^p(\Sigma_{\alpha,\beta})$ is dense in $\mathcal{M}_{\sigma}(\Sigma_{\alpha,\beta})$.
\end{lem}

\subsection{Hofbuaer's Markov diagram for ($\alpha$-$\beta$)-shifts}

The Hofbauer's Markov diagram is a countable oriented graph, whose vertices are subsets of $\Sigma^+_{\alpha,\beta}$.
Let $D\subset\Sigma_{\alpha,\beta}^+$ be a closed subset with $D\subset [i]$ for some $1\le i\le k$.
We say that a non-empty closed subset $C\subset\Sigma_{\alpha,\beta}^+$ is a \textit{successor} of $D$ if
$C=[j]\cap\sigma(D)$ for some $1\le j\le k$.

Now we define a set $\mathcal{D}=\mathcal{D}_{\alpha,\beta}$ of vertices by induction. First, we set
$\mathcal{D}_0:=\{[1],\ldots,[k]\}$. If $\mathcal{D}_n$ is defined for $n\ge 0$, then we set
$$\mathcal{D}_{n+1}:=\mathcal{D}_n\cup \{C:\text{there exists }D\in\mathcal{D}_n\text{ so that }C\text{ is a
successor of }D\}.$$
We note that $\mathcal{D}_n$ is a finite set for each $n$ since the number of successors of $D$ is
at most $k$ by the definition. Finally, we set
$$\mathcal{D}:=\bigcup_{n\ge 0}\mathcal{D}_n.$$
We insert an arrow from every
$D\in\mathcal{D}$ to all its successors. We write $D\rightarrow C$ if
$C$ is a successor of $D$.
We call the countable oriented graph $(\mathcal{D},\rightarrow)$ the \textit{Hofbauer's Markov Diagram}.
Let $\mathcal{C}\subset \mathcal{D}$.
We say that a sequence $C_1\cdots C_n$ is a \textit{finite path} in $\mathcal{C}$ if
$C_1,\ldots,C_n\in\mathcal{C}$ and $C_i\rightarrow C_{i+1}$ for $1\le i\le n-1$.
A two-sided infinite path in $\mathcal{C}$ is defined analogously.
We denote by $\mathcal{X}(\mathcal{C})$ the set of all two-sided infinite paths in $\mathcal{C}$.
$\mathcal{C}$ is said to be \textit{irreducible} if for any $C,D\in\mathcal{C}$,
there is a finite path $C_1C_2\cdots C_n$ in $\mathcal{C}$ with $C_1=C$ and $C_n=D$.
%and if every subset of $\mathcal{D}$, which contains $\mathcal{C}$ strictly, does not have this property.
Note that $\mathcal{X}(\mathcal{C})$ is a countable two-sided Markov shift for any subset $\mathcal{C}\subset\mathcal{D}$.
Indeed, if we define a matrix $M=(M_{D,C})_{(D,C)\in\mathcal{C}^2}$ by
$$M_{D,C}=
\left\{
\begin{array}{ll}
1 & (D\rightarrow C), \\
0 & (\text{otherwise}).
\end{array}
\right.$$
then we have $\mathcal{X}(\mathcal{C})=\Sigma_M$.
Moreover, $\mathcal{X}(\mathcal{C})$ is transitive if $\mathcal{C}$ is irreducible.
Since $\beta>2$, we can easily see that $\mathcal{D}=\mathcal{D}_{\alpha,\beta}$ is irreducible.
Hence, we have the following lemma:

\begin{lem}
\label{markov}
$\mathcal{X}(\mathcal{D})$ is a transitive countable Markov shift.
\end{lem}

The irreducibility of $\mathcal{D}$ also implies the following:

\begin{lem}
\label{connecting-time}
Let $\mathcal{F},\mathcal{F}'\subset\mathcal{D}$ be two finite subsets.
Then there is an integer $t=t(\mathcal{F},\mathcal{F}')\ge 1$ such that
for any $C\in\mathcal{F}$ and $D\in\mathcal{F}'$, we can find a finite path
$C_0C_1\cdots C_tC_{t+1}$ in $\mathcal{D}$ such that $C_0=C$ and $C_{t+1}=D$.
\begin{proof}
For $C\in\mathcal{F}$ and $D\in\mathcal{F}'$, we set
$$t_1(C):=\min\{n\ge 1:C_1\cdots C_n\text{ is a finite path in $\mathcal{D}$ with }C_1=C\text{ and }C_n=[2]\},$$
$$t_2(D):=\min\{n\ge 1:C_1\cdots C_n\text{ is a finite path in $\mathcal{D}$ with }C_1=[2]\text{ and }C_n=D\}.$$
Since $\mathcal{D}$ is irreducible, we have $t_1(C),t_2(D)<\infty$.
Then, we define
$t(\mathcal{F},\mathcal{F}'):=\max\{t_1(C)+t_2(D):C\in\mathcal{F},D\in\mathcal{F}'\}$, which is finite
since $\mathcal{F}$ and $\mathcal{F}'$ are finite.

Take any $C\in\mathcal{F}$ and $D\in\mathcal{F}'$.
For the notational simplicity, we set
$t:=t(\mathcal{F},\mathcal{F}')$, $t_1:=t_1(C)$ and $t_2:=t_2(D)$.
By the definition of $t_1$, there exists a finite path $C_0C_1\cdots C_{t_1-1}$ in $\mathcal{D}$ such that
$C_0=C$ and $C_{t_1-1}=[2]$. Similarly, we can find a finite path $C_{t-t_2+2}\cdots C_{t+1}$ in $\mathcal{D}$
such that $C_{t-t_2+2}=[2]$ and $C_{t+1}=D$. For $t_1\le j\le t-t_2+1$, we set $C_j:=[2]$.
Since $\beta>2$, it is clear that $[2]\rightarrow [2]$.
Hence $C_0\cdots C_{t+1}$ is a finite path in $\mathcal{D}$, which proves the lemma.
\end{proof}
\end{lem}

It is known that there is a deep connection between two subshifts ($\mathcal{X}(\mathcal{D}),\sigma)$
and $(\Sigma_{\alpha,\beta},\sigma)$.
To explain this connection, we define $\Psi\colon \mathcal{X}(\mathcal{D})\to \{1,\ldots,k\}^{\mathbb{Z}}$ as follows.
For each $D\in\mathcal{D}$, we can find a unique $1\le i\le k$ so that $D\subset [i]$ holds.
Hence, we can define $\psi(D)=i$ if $D\subset [i]$.
For $(D_i)_{i\in\mathbb{Z}}\in \mathcal{X}(\mathcal{D})$, we set
$$\Psi((D_i)_{i\in\mathbb{Z}}):=(\psi(D_i))_{i\in\mathbb{Z}}.$$
The following lemma, which is important to show Theorem \ref{main} is proved by Hofbauer
in \cite{H2} (see also \cite[Appendix]{Bu}).

\begin{lem}(\cite[Lemmas 2 and 3]{H2})
\label{hofbauer-extension}
The map $\Psi\colon\mathcal{X}(\mathcal{D})\to \Sigma_{\alpha,\beta}$ is continuous
and $\sigma\circ\Psi=\Psi\circ\sigma$. Moreover, we can find
two $\sigma$-invariant subsets $\tilde{N}\subset\mathcal{X}(\mathcal{D})$ and $N\subset\Sigma_{\alpha,\beta}$ satisfying the following properties:
\begin{itemize}
\item
The restriction map $\Psi\colon\mathcal{X}(\mathcal{D})\setminus \tilde{N}\to \Sigma_{\alpha,\beta}\setminus N$
is bijective and bi-measurable. 

\item
$N$ has no periodic point.

\item
For any invariant measure $\mu\in\mathcal{M}_{\sigma}(\Sigma_{\alpha,\beta})$ with $\mu(N)=1$, we have
$h(\mu)=0$.

\item
For any invariant measure $\tilde{\mu}\in\mathcal{M}_{\sigma}(\mathcal{X}(\mathcal{D}))$ with
$\tilde{\mu}(\tilde{N})=1$, we have $h(\tilde{\mu})=0$.
\end{itemize}
\end{lem}

\section{Proof of Theorem \ref{main}}

In this section, we give a proof of Theorem \ref{main}. Take an arbitrary
$\mu\in\mathcal{M}_{\sigma}(\Sigma^+_{\alpha,\beta})$. Without loss of generality we may assume $h(\mu)>0$.
We begin with the following proposition:

\begin{prop}
\label{entropy-dense}
For any $\epsilon>0$ and any neighborhood $\mathcal{U}$ of $\mu$
in $\mathcal{M}(\Sigma_{\alpha,\beta})$, there exist a finite set $\mathcal{F}\subset \mathcal{D}$ and
an ergodic measure $\rho\in\mathcal{M}_{\sigma}^e(\pi(\Psi(\mathcal{X}(\mathcal{F}))))$ such that
$\rho\in \mathcal{U}$ and $|h(\mu)-h(\rho)|\le\epsilon$.
Here $\pi\colon\Sigma_{\alpha,\beta}\to\Sigma_{\alpha,\beta}^+$;
$(x_i)_{i\in\mathbb{Z}}\mapsto (x_i)_{i\in\mathbb{N}}$ be a canonical projection and
$\Psi\colon\mathcal{X}(\mathcal{D})\to\{1,\ldots,k\}^{\mathbb{Z}}$ is as in \S2.
\begin{proof}
First, we define two maps $\pi_{\ast}\colon \mathcal{M}_{\sigma}(\Sigma_{\alpha,\beta})
\to\mathcal{M}_{\sigma}(\Sigma_{\alpha,\beta}^+)$ and\\
$\Psi_{\ast}\colon\mathcal{M}_{\sigma}(\mathcal{X}(\mathcal{D}))\to
\mathcal{M}_{\sigma}(\Sigma_{\alpha,\beta})$
by $\pi_{\ast}(\mu^-):=\mu^-\circ\pi^{-1}$ and $\Psi_{\ast}(\tilde{\mu}):=\tilde{\mu}\circ\Psi^{-1}$.
Then it is well-known that $\pi_{\ast}$ is a homeomorphism and
$h(\mu^-)=h(\pi_{\ast}(\mu^-))$ for any $\mu^-\in\mathcal{M}_{\sigma}(\Sigma_{\alpha,\beta})$.
Moreover, it follows from Lemma \ref{hofbauer-extension} that
$\Psi_{\ast}$
is continuous, and the restriction map $\Psi_{\ast}\colon\mathcal{M}_{\sigma}(\mathcal{X}(\mathcal{D})\setminus \tilde{N})
\to \mathcal{M}_{\sigma}(\Sigma_{\alpha,\beta}\setminus N)$ is bijective, bi-measurable, and
$h(\tilde{\mu})=h(\Psi_{\ast}(\tilde{\mu}))$ holds for any $\tilde{\mu}\in
\mathcal{M}_{\sigma}(\mathcal{X}(\mathcal{D})\setminus\tilde{N})$.
%$$\Psi_{\ast}(\mathcal{M}_{\sigma}(\mathcal{X}(\mathcal{D})))=\{\mu\in
%\mathcal{M}_{\sigma}(\Sigma_{\alpha,\beta}):\mu(\Sigma_{\alpha,\beta}\setminus N)=1\}.$$
Here $N\subset \Sigma_{\alpha,\beta}$ and $\tilde{N}\subset \mathcal{X}(\mathcal{D})$ are as in Lemma \ref{hofbauer-extension},
$$\mathcal{M}_{\sigma}(\mathcal{X}(\mathcal{D})\setminus\tilde{N}):=\{\tilde{\mu}\in\mathcal{M}_{\sigma}(\mathcal{X}(\mathcal{D})):\tilde{\mu}(\tilde{N})=0\}\text{ and}$$
$$\mathcal{M}_{\sigma}(\Sigma_{\alpha,\beta}\setminus N):=\{\mu^-\in\mathcal{M}_{\sigma}(\Sigma_{\alpha,\beta}):
\mu^-(N)=0\}.$$
We set $\mu^-:=\pi_{\ast}^{-1}(\mu)\in\pi_{\ast}^{-1}(\mathcal{U})$.
Then we can find $0\le a\le 1$, $\mu^-_1,\mu_2^-\in\mathcal{M}_{\sigma}(\Sigma_{\alpha,\beta})$
such that
$\mu^-=a\mu^-_1+(1-a)\mu^-_2$, $\mu_1^-(\Sigma_{\alpha,\beta}\setminus N)=1$ and $\mu_2^-(N)=1$.
We note that $h(\mu^-_2)=0$ by Lemma \ref{hofbauer-extension}.
By Lemma \ref{per-dense}, we can find $\nu^-_2\in\mathcal{M}^p(\Sigma_{\alpha,\beta})$
such that $\nu^-:=a\mu^-_1+(1-a)\nu^-_2\in\pi_{\ast}^{-1}(\mathcal{U})$.
It is clear that $h(\nu^-)=h(\mu^-)=h(\mu)$.
Since $N$ has no periodic point, $\nu^-_2(\Sigma_{\alpha,\beta}\setminus N)=1$, which also implies that
$\nu^-(\Sigma_{\alpha,\beta}\setminus N)=1$.
Therefore, we can define $\tilde{\nu}:=\Psi_{\ast}^{-1}(\nu^-)$.
Then we have $\tilde{\nu}\in U$, where we set
$U:=\Psi_{\ast}^{-1}((\pi_{\ast})^{-1}(\mathcal{U}))$, which is open since $\pi_{\ast}\circ\Psi_{\ast}$
is continuous.
Note that $\tilde{\nu}$ is supported on $\mathcal{X}(\mathcal{D})$ and
$\mathcal{X}(\mathcal{D})$ is a transitive two-sided
Markov shift by Lemma \ref{markov}.
Thus, it follows from Lemma \ref{takahasi} that we can find
a finite set $\mathcal{F}\subset \mathcal{D}$ and an ergodic measure
$\tilde{\rho}\in U$ supported on $\mathcal{X}(\mathcal{F})$
such that $|h(\tilde{\nu})-h(\tilde{\rho})|\le\epsilon$.
Without loss of generality, we may assume that $h(\tilde{\rho})>0$, which implies that
$\tilde{\rho}(\tilde{N})=0$.
We define an ergodic measure $\rho$ on $\Sigma_{\alpha,\beta}^+$ by $\rho:=\pi_{\ast}(\Phi_{\ast}(\tilde{\rho}))$.
Clearly, we have $\rho\in \mathcal{M}_{\sigma}^e(\pi(\Psi(\mathcal{X}(\mathcal{F}))))\cap \mathcal{U}$.
Moreover, since $\tilde{\rho}\in\mathcal{M}_{\sigma}(\mathcal{M}(\mathcal{X}(\mathcal{D})\setminus\tilde{N}))$,
we have $h(\tilde{\rho})=h(\rho)$.
Hence we have $|h(\mu)-h(\rho)|\le \epsilon$,
which proves the proposition.
\end{proof}
\end{prop}

Let $\epsilon>0$ and choose a strictly decreasing sequence $\{\epsilon_k\}_{k\ge 1}$ of positive real
numbers so that $\epsilon_k\rightarrow 0$ and $\epsilon_1\le \epsilon$.
Take a compatible metric $D$ on $\mathcal{M}(\Sigma_{\alpha,\beta}^+)$ so that
$D(\nu,\nu')\le 1$ holds for any $\nu,\nu'\in\mathcal{M}(\Sigma_{\alpha,\beta}^+)$.
Then it follows from Proposition \ref{entropy-dense} that
for any $k\ge 1$, we can find a finite $\mathcal{F}_k\subset \mathcal{D}$ and
$\mu_k\in\mathcal{M}_{\sigma}^e(\pi(\Psi(\mathcal{X}(\mathcal{F}_k))))$ such that
$D(\mu,\mu_k)\le\epsilon_k$ and $h(\mu)-\epsilon\le h(\mu_k)$.
We set $X_k:=\pi(\Psi(\mathcal{X}(\mathcal{F}_k)))$.
Since $\mu_k$ is an ergodic measure on $X_k$, it follows from \cite[Propositions 2.1 and 4.1]{PS} that
we can find $l_k\in\mathbb{N}$ so that
$$\#\left\{w\in\mathcal{L}_{l_k}(X_k):
D\left(\mu_k,\frac{1}{l_k}\sum_{j=0}^{l_k-1}\delta_{\sigma^j(x)}\right)\le\epsilon_k,x\in[w]\right\}\ge
e^{l_k(h(\mu_k)-\epsilon)}$$
$$\text{and }\frac{\max\{t(\mathcal{F}_k,\mathcal{F}_k),t(\mathcal{F}_k,\mathcal{F}_{k+1})\}}{l_k}\le\epsilon_k$$
hold. Here $t(\mathcal{F},\mathcal{F}')$ is as in Lemma \ref{connecting-time}.
For the notational simplicity, we set
$$\Gamma_k:=\left\{w\in\mathcal{L}_{l_k}(X_k):
D\left(\mu_k,\frac{1}{l_k}\sum_{j=0}^{l_k-1}\delta_{\sigma^j(x)}\right)\le\epsilon_k,x\in[w]\right\}.$$
We also choose sequence $\{L_k\}$ such that
$$\max\left\{l_{k+1},\sum_{j=1}^{k-1}l_jL_j\right\}\le\epsilon_k\sum_{j=1}^kl_jL_j$$
holds.
For $k\ge 1$ with
$$k=L_1+\cdots +L_{j-1}+q$$
for some $j\ge 1$ and $1\le q\le L_j$, we define four sequences
$\{l_k'\}$, $\{\mathcal{F}_k'\}$, $\{n_k\}$ and $\{\Gamma_k'\}$ as
$$l_k':=l_j,\mathcal{F}_k':=\mathcal{F}_j,n_k:=l_k'+t(\mathcal{F}_k',\mathcal{F}_{k+1}'),\Gamma_k':=\Gamma_j$$
and define
$$G:=\bigcap_{k=1}^{\infty}\bigcup_{w\in\Gamma_k'}\sigma^{-(n_1+\cdots+n_{k-1})}[w].$$
Then by the definition, we can easily to see the following properties:
\begin{itemize}
\item
$\displaystyle \lim_{k\rightarrow\infty}\frac{n_1+\cdots n_k}{n_1+\cdots n_{k-1}}=1$.

\item
$\displaystyle\lim_{j\rightarrow\infty}\frac{n_1+\cdots +n_{m_j}}{n_1+\cdots+ n_{m_{j+1}}}=0$,
where we set $m_j:=L_1+\cdots +L_j$.

\item
For any $k\ge 1$, $w\in\Gamma_k'$, and $x\in [w]$, we have
$$D\left(\mu,\frac{1}{n_k}\sum_{j=0}^{n_k-1}\delta_{\sigma^j(x)}\right)\le 3\epsilon_k.$$

\item
For any $k\ge 1$, we have $\#\Gamma'_k\ge\exp\{n_k(\frac{1}{1+\epsilon}(h(\mu)-2\epsilon))\}$.

\item
For any $\displaystyle (w^1w^2\cdots)\in\prod_{k=1}^{\infty}\Gamma_k'$,
$\displaystyle G(w^1w^2\cdots):=\bigcap_{k=1}^{\infty}\sigma^{-(n_1+\cdots +n_{k-1})}[w]$
is a non-empty closed set. Moreover, we can write
$$G=\bigcup\left\{G(w^1w^2\cdots):(w^1w^2\cdots)\in\prod_{k=1}^{\infty}\Gamma_k'\right\}.$$

\end{itemize}

These properties enable us to show $G\subset G_{\mu}$ and $h_{\rm top}(\sigma,G)\ge
\frac{1}{1+\epsilon}(h(\mu)-2\epsilon)$
in a similar way to the proof of \cite[Lemma 5.1]{PS2}.
These together with Lemma \ref{upper} imply Theorem \ref{main}.

\vspace*{3mm}

\noindent
\textbf{Acknowledgement.}~ 
The author was partially supported by JSPS KAKENHI Grant Number 18K03359.

\end{document}